\documentclass{amsart}
\usepackage{amssymb}
\title{Norm closure of classical pseudodifferential operators does not contain H\"ormander's class}
\author{Severino T. Melo}
\date{}
\newtheorem{thm}{Theorem}

\newtheorem{lem}{Lemma}
\newtheorem{cor}{Corollary}

\begin{document}
%
%
%
%
\newcommand{\op}{operator}
\newcommand{\ops}{operators}
\newcommand{\psd}{pseudo\-dif\-fer\-en\-tial}
\newcommand{\Psd}{Pseudo\-dif\-fer\-en\-tial}
\newcommand{\pf}{{\em Proof}: }
\newcommand{\ifoi}{if and only if}
\newcommand{\cst}{C$^{*}$}
%
%
%
%
\newcommand{\m}{{\mathbb M}}
\newcommand{\z}{{\mathbb Z}}
\newcommand{\co}{{\mathbb C}}
\newcommand{\re}{{\mathbb R}}
\newcommand{\so}{{\mathbb S}^1}
\newcommand{\sn}{{\mathbb S}^n}
\newcommand{\rn}{{\mathbb R}^{n}}
\newcommand{\sx}{S^{*}X}
\newcommand{\tix}{T^*X}
\newcommand{\bix}{B^*X}
%
%
%
%
\newcommand{\pscl}{\Psi_{cl}^{0}(X)}
\newcommand{\psh}{\Psi^{0}(X)}
\newcommand{\cix}{C^{\infty}(X)}
\newcommand{\cisx}{C^{\infty}(\sx)}
\newcommand{\citx}{C^{\infty}(\tix)}
\newcommand{\cis}{C^{\infty}(\so)}
\newcommand{\cibx}{C^{\infty}(\bx)}
\newcommand{\cisbx}{C^{\infty}(\sbx)}
\newcommand{\h}{{\mathfrak H}}
\newcommand{\lh}{{\mathfrak L}(\h)}
\newcommand{\as}{{\mathcal A}}
\newcommand{\ac}{{\mathfrak A}}
\newcommand{\Sc}{{\mathfrak S}}
\newcommand{\Cc}{{\mathfrak C}}
\newcommand{\bc}{{\mathfrak B}}
\newcommand{\acc}{{\mathfrak A}_{cl}}
\newcommand{\kc}{{\mathfrak K}}
%
%
%
%
\newcommand{\sig}{\bar\sigma}
\newcommand{\cqd}{\hfill$\Box$}
\newcommand{\supp}{{\rm supp}}
\newcommand{\du}{d\!u}
\newcommand{\dv}{d\!v}
\newcommand{\dx}{d\!x}
\newcommand{\dy}{d\!y}
\newcommand{\dS}{d\!S}
\newcommand{\dxi}{d\!\xi}

\begin{abstract}
The $L^2$-norm closure $\acc$\ of the set $\pscl$\ of all zero-order classical \psd\ 
\ops\ on a closed manifold $X$\ is proven not to contain H\"ormander's class $\psh$. 
A set of generators for $\acc$\ smaller than all of $\pscl$\ is also described. 
\end{abstract}
\maketitle
\section{Introduction}

Let $X$\ denote a closed manifold of dimension $n$. As in his book \cite{H}, we denote by 
$\psh$\ the set of all zero-order elements in H\"ormander's class of \psd\ \ops\ (of type 
$\rho=1$\ and $\delta=0$). The subclass of classical \psd\ \ops\ of order zero, i.e., those 
whose local symbols possess asymptotic expansions in homogeneous components of decreasing 
degree, will be here denoted by $\pscl$\ (this is denoted by $\Psi_{phg}^{0}(X)$\ in \cite{H}).

The Hilbert space of square integrable functions over $X$\ (equipped with a measure defined 
by a positive smooth density) will be denoted by $\h$, and its algebra of bounded \ops\ by 
$\lh$. The classes of \psd\ \ops\ $\psh$\ and $\pscl$\ are adjoint-invariant subalgebras of 
$\lh$. Their norm closures will be denoted by $\ac$\ and $\acc$, respectively. Both contain the 
ideal $\kc$\ of compact \ops, and their quotients by it are unital commutative \cst-algebras; 
hence isomorphic to spaces of continuous functions over compact Hausdorff spaces. That $\ac$\ 
is not contained in $\acc$\ follows from the fact that these spaces, the Gelfand spaces of 
$\ac/\kc$\ and of $\acc/\kc$, are not the same; as we prove in Sections \ref{gs}\ and 
\ref{gsc}. 

Our proofs strongly depend on H\"ormander's estimate for the norm, modulo compact operators,
of a \psd\ \op\ (\cite{Hp}, Theorem~3.3). For classical \psd\ \ops, that estimate had  
already appeared in \cite{KN}, Theorem~A.4; and had been proven, for singular integral
\ops, by Gohberg \cite{G}\ and Seeley \cite{S}.

In Section \ref{ca}, we show that $\acc$\ is actually a {\em comparison algebra}, as defined by 
Cordes \cite{C2}. This means that it is equal to the smallest \cst-subalgebra of $\lh$\ 
containing all \ops\ of multiplication by smooth functions and all operators of the form
$L(1-\Delta)^{-1/2}$, where $L$\ is any first-order linear differential \op\ with smooth
coefficients and $\Delta$\ is a certain second-order non-positive elliptic differential operator 
on $X$. That can be proven again by a Gelfand space argument.

On page 508 in the first of a series of papers about their index theorem, Atiyah and
Singer \cite{AS}\ state: ``The term \psd\ \op\ is applied, in different places, to slightly 
different classes of \ops. For our purposes any one of these classes would be equally good.
In fact, we shall eventually form a closure of this class and, by that stage, {\em any 
differences would disappear}. Perhaps the largest and most natural class is that given 
by H\"ormander \cite{Hp} (...)''. The norm closure of the set of all \psd\ \ops\ considered by 
Atiyah and Singer in \cite{AS}\ (the subclass of H\"ormander's class defined on page 509) is 
however smaller than the closure of the class introduced in \cite{Hp}; whose zero-order
subclass was later denoted $\psh$. Indeed, the results in Section~\ref{gsc}\ show, for example, 
that the image of the continuous extension of the principal symbol to the closure of $\psh$\ 
cannot be identified with the space of continuous functions on the co-sphere bundle, as 
described in the Remark right after (5.2) in \cite{AS}.

Taylor proved \cite{T}\ that a bounded operator $A$\ on $L^2(\sn)$\ belongs to $\Psi^0(\sn)$\ 
\ifoi\ the mapping 
\[
SO_e(n+1,1)\ni g\mapsto U(g)AU(g)^{-1}\in{\mathfrak L}(L^2(\sn))
\]
is smooth, where $U$\ denotes the unitary representation of the conformal group $SO_e(n+1,1)$\ 
on $L^2(\sn)$\ induced by its natural action on the sphere $\sn$\ (for simplicity, only this 
particular case of his result is stated here). The fact that the norm closure of the classical 
\psd\ \ops\ is smaller than that of H\"ormander's class makes it plausible to expect that, 
perhaps, the classical \psd\ \ops\ can also be characterized by a similar smoothness condition, 
for a larger group. On the other direction, the class of smooth \ops\ for the action of (the 
subgroup of the conformal group) $SO(n)$\ contains \ops\ without the pseudolocal property 
\cite{CM,M}. We refer to \cite{P}\ for possible applications of these smoothness 
characterizations of \psd\ \ops, the first of which was given by Cordes \cite{Cs}. 
    
\section{The Gelfand space of $\ac/\kc$}\label{gs}

The ``basic calculus'' (see \cite{H}, Section 18.1; or \cite{Tr}, Sections I.3 and I.4) of 
\psd\ \ops\ implies that the commutator of any two elements of $\psh$\ is of order $-1$. 
Operators of order $-1$\ are bounded from $\h$\ into the Sobolev space $H^1(X)$, which 
compactly embeds in $\h$. Since $\psh$\ is dense in $\ac$, all commutators of $\ac$\ then 
belong to the compact ideal $\kc$.

Any integral \op\ with smooth kernel is a \psd\ \op\ of order $-\infty$\ on $X$. In particular, 
any finite rank \op\ with image and kernel contained in $\cix$\ belongs to $\psh$. A simple 
density argument gives that all finite rank \ops\ belong to $\ac$, and then 
$\kc\subset\ac$. The quotient $\ac/\kc$\ is therefore commutative. It is, moreover, a 
\cst-algebra (\cite{Mu}, 3.1.4). 

Any unital commutative \cst-algebra $\bc$\ is isomorphic to the algebra of continuous 
functions on the set $\m$\ of its nonzero multiplicative linear functionals (which are 
automatically continuous of norm $1$) equipped with the weak-* topology. This is Gelfand's 
theorem (\cite{Mu}, 2.1.10), and we will call $\m$\ the Gelfand space of $\bc$. That theorem 
also states that the {\em Gelfand mapping}\ for $\bc$,
\[
\begin{array}{rcl}
\bc&\longrightarrow&C(\m)\\
x&\longmapsto&(\m\ni\varphi\mapsto \varphi(x))
\end{array},
\]
is an isomorphism. 

In this Section, we show that the Gelfand space of $\ac/\kc$\ are
the points at infinity of the compactification of the cotangent bundle determined
by the zero-order symbols, and that the Gelfand mapping then corresponds to the restriction 
of the principal symbol over those points.

Let $\Sc$\ denote the closure in sup-norm of $S^0(\tix)$, the set of symbols of order 
zero on the cotangent bundle $\tix$\ (as defined in \cite{H}, page 85). This is a 
\cst-subalgebra of the set of all continuous bounded functions on $\tix$. Let $\m_s$\ denote the  
Gelfand space of $\Sc$. Evaluation at each point defines a mapping from $\tix$\ to $\m_s$. It 
follows from the definition of weak-* topology that this mapping is continuous. It follows from 
Urysohn's Lemma that it is a homeomorphism onto its dense image ($\m_s$\ is a compact Hausdorff 
space). In other words, $\m_s$\ is a compactification of $\tix$. Regarding this embedding of 
$\tix$\ into $\m_s$\ as an inclusion, let us define  
\[
\m_h=\m_s\setminus\tix.  
\]
The unique continuous extension of any $a\in S^0(\tix)$\ to $\m_s$\ will also be
denoted by $a$.

\begin{lem}\label{smo}
Any symbol of order $-1$, $a\in S^{-1}(\tix)$, vanishes on $\m_h$.
\end{lem}

\pf 
Given $m\in \m_h$, there is a net $v_\alpha\in\tix$\ converging to $m$. Denoting by 
$\pi:\tix\to X$\ the bundle projection and, if necessary, passing to a subnet, we may suppose 
that every $x_\alpha=\pi(v_\alpha)\in X$\ belongs to the domain of a chart 
$\chi:U\to\tilde U$\ and that $x_\alpha$\ converges to an $x_0\in U$\ (we have used that
any net in the compact space $X$\ has a convergent subnet; see \cite{RS}, Theorem IV.3). 

Let $\phi:T^*U\to\tilde U\times\rn$\ denote the bundle trivialization induced by $\chi$, and 
let $\phi(v_\alpha)=(y_\alpha,\xi_\alpha)$. I claim that $|\xi_\alpha|\to\infty$. Indeed, if 
that were not true, $\xi_\alpha$\ would have a subnet $\xi_{F(\alpha)}$\ converging to some 
$\xi_0\in\rn$. Then $v_{F(\alpha)}$\ would converge to $\phi^{-1}(\chi(x_0),\xi_0)$. 
That would contradict the assumption that $\lim v_\alpha$\ is a point in 
$\m_h=\m_s\setminus\tix$.

The pushforward of $a$\ by $\phi$\ is, by definition, a local symbol
of order $-1$, $\phi_*a\in S^{-1}(\tilde U\times\rn)$. It also no loss of generality to 
suppose that the net with which we started is such that there is a compact subset of $\tilde U$\
containing every $y_\alpha$; and then 
$\phi_*a(y_\alpha,\xi_\alpha)=a(\phi^{-1}(y_\alpha,\xi_\alpha))$\ converges to zero. This 
proves that $a(m)=\lim a(v_\alpha)=0$. \cqd

\begin{lem}\label{topo} For all $a\in S^0(\tix)$\ we have:
\begin{equation}
\label{limsup}
\sup\{|a(m)|;\,m\in\m_h\}\,=\,\lim_{R\to\infty}\sup\{|a(v)|;\,v\in\tix,\,|v|^\prime>R\},
\end{equation}
where $|v|^\prime$\ denotes the evaluation at $v$\ of the norm defined on $T_{\pi(v)}^*X$\ by a 
Riemannian metric on $X$. 
\end{lem}

\pf Since $\m_h$\ is compact, there is an $m_0\in\m_h$\ such that $a(m_0)$\ equals the
left-hand side of (\ref{limsup}), which we denote by $||a|_{\m_{h}}||_\infty$. Let $v_\alpha$\
denote a net in $\tix$\ which converges to $m_0$. We may suppose, as in the proof of 
Lemma~\ref{smo}, that $v_\alpha$\ is locally given by $(y_\alpha,\xi_\alpha)$, with
$|\xi_\alpha|\to\infty$, $|\cdot|$\ denoting the euclidean norm on $\rn$. It is 
easy to see that $|v_\alpha|^\prime\to\infty$. If we define $R_\alpha=|v_\alpha|^\prime/2$,
we then have $\lim_{\alpha}R_\alpha=\infty$. Since 
\[
|a(v_\alpha)|\,\leq\,\sup\{|a(v)|;\,|v|^\prime>R_\alpha\,\},
\] 
we have:
\[
||a|_{\m_{h}}||_\infty\,=\,\lim_{\alpha}|a(v_\alpha)|\,\leq
\]
\[
\lim_{\alpha}\,\sup\{|a(v)|;|v|^\prime>R_\alpha\}\,=\,
\lim_{R\to\infty}\sup\{|a(v)|;|v|^\prime>R\}.
\]

Suppose the above inequality is strict. Then there is $\delta>0$\ and a sequence $v_k\in\tix$\
such that $|a(v_k)|$\ converges to $||a|_{\m_{h}}||+\delta$\ and $|v_k|^\prime\to\infty$. Let 
$v_\alpha$\ be a convergent subnet of that sequence. Its limit $m_0$\ cannot belong to $\tix$, 
since $|v_\alpha|^\prime\to\infty$. That gives $|a(m_0)|>||a|_{\m_{h}}||_\infty$\ for
some $m_0\in\m_h$, which is a contradiction. \cqd

The principal symbol of an operator in $\psh$\ is defined modulo negative-order terms (\cite{H}, 
pages 86 and 87). That gives  a surjective homomorphism from $\psh$\ onto 
$S^0(\tix)/S^{-1}(\tix)$. By Lemma~\ref{smo}, the restrictions to $\m_h$\ of any two
representatives of the principal symbol of a given $A\in\psh$\ coincide. This defines a 
mapping 
\[
\sigma:\psh\longrightarrow C(\m_h),
\]
which is an algebra homomorphism (\cite{H}, Theorem 18.1.23). It also follows from the 
basic calculus that $\sigma(A^*)=\overline{\sigma(A)}$, for all $A\in\psh$, with
$A^*$\ denoting the Hilbert-space adjoint of $A$, and bar denoting complex conjugate.

\begin{thm}\label{gh}
The *-homomorphism $\sigma$\ is bounded. Its continuous extension to $\ac$, still denoted
by $\sigma$, is surjective and has kernel equal to the compact ideal $\kc$.
\end{thm} 
\pf 
The proof of \cite{Hp}, Theorem 3.3 (see the Remark at the end of that proof), shows that,
for every $A\in\psh$,
\begin{equation}
\label{estimate}
\inf\{||A+K||;\,K\in\kc\}\,=\,\lim_{R\to\infty}\sup\{|a(v)|;\,v\in\tix,\,|v|^\prime>R\},
\end{equation}
which equals $||a|_{\m_{h}}||_\infty$, by Lemma~\ref{topo}. 

Since the left-hand side of (\ref{estimate}) is bounded by $||A||$, we obtain that $\sigma$\ 
extends to a bounded \op\ between the Banach spaces $\ac$\ and $C(\m_h)$. It also follows from 
(\ref{estimate}) that the kernel of that extension is equal to $\kc$, and that 
the mapping induced on the quotient is an isometry; and hence has closed image. By 
definition, $\{\sigma(A);A\in\psh\}$\ is dense in $C(\m_h)$, and hence $\sigma$\ is onto.
\cqd

\begin{cor}\label{cgh}The mapping
\[
\begin{array}{rcl}
\ac/\kc&\longrightarrow&C(\m_h)\\
\,[A] &\longmapsto&\sigma(A)
\end{array}
\]
is a \cst-algebra isomorphism {\em (}$[A]$\ denotes the class of $A\in\ac$\
in the quotient $\ac/\kc${\em )}. 
\end{cor}

\section{The Gelfand space of $\acc/\kc$}\label{gsc}

An operator $A\in\psh$\ is {\em classical}, or {\em polihomogeneous}, if for every chart
$\chi:U\to\tilde U$, the pushforward $\chi_*A$,
\[
C_0^\infty(\tilde U)\ni u\mapsto (\chi_*A)u=[A(u\circ\chi)]\circ\chi^{-1}\in C^\infty(\tilde U),
\]
is equal, modulo an integral \op\ with smooth kernel, to $\tilde a(x,D)$, for some 
$\tilde a\in S^0(\tilde U\times\rn)$\ possessing an asymptotic expansion
\[
\tilde a(x,\xi)\,\sim\,\sum_{j=0}^{\infty}a_{j}(x,\xi),
\]
where each $a_j$\ satisfies $a_j(x,t\xi)=t^{-j}a_j(x,\xi)$, for every $t>0$\ and every $\xi$\ 
with $|\xi|>1$. Although $\tilde a$\ and $a_j$\ for $j>0$\ are only locally defined, $a_0$\
is a globally defined function on the cotangent bundle (see, for example, the comments at
the end of Section I.5 in \cite{Tr}). That gives a representative of the principal 
symbol of $A$\ satisfying, for some $R_0>0$, 
\begin{equation}
\label{homogeneous}
a(tv)=a(v),\ \mbox{if}\ t>0\ \mbox{and}\ |v|^\prime>R_0. 
\end{equation}
Moreover, it is possible to arrange $R_0$\ to be smaller than one. 
For such an $a$, $\sup\{|a(v)|;|v|^\prime>R\}$\ is independent of $R$, for $R\geq 1$. 
If $a_1$\ is another principal symbol of $A$\ satisfying (\ref{homogeneous}) with $R_0<1$,
then $a_1$\ and $a$\ coincide on $\{|v|^\prime\geq 1\}$, by Lemmas~\ref{smo}\ and \ref{topo}.
The following theorem then follows from the estimate (\ref{estimate}).

\begin{thm}
If $A\in\psh$\ is classical (denoted $A\in\pscl$) and $a\in\citx$\ is a principal 
symbol of $A$\ satisfying $a(tv)=a(v)$\ for all $t>0$\ and for all $v\in\tix$\ such 
that $|v|^\prime\geq 1$, then
\begin{equation}
\label{cosphere}
\inf\{||A+K||;\,K\in\kc\}\,=\,\sup\{|a(v)|,\,v\in\sx\},
\end{equation}
where $\sx$\ denotes the co-sphere bundle, $\sx=\{v\in\tix;|v|^\prime =1\}$.
\label{gc}\end{thm}

\begin{cor}\label{cgc} The mapping $[A]\mapsto a|_{\sx}$, for $A\in\pscl$\ and
$a$\ as in the statement of Theorem~\ref{gc}, extends to a \cst-algebra isomorphism 
$\hat\sigma:\acc/\kc\to C(\sx)$.
\end{cor}
\pf That this defines an isometric \cst-algebra homomorphism follows exactly as in the proof
of Theorem~\ref{gc}. To prove that it is surjective, it is enough to show that its image
contains the dense subspace $\cisx$. Given $b\in\cisx$, let $a\in\citx$\ be such that
(\ref{homogeneous}) holds for some $R_0<1$. The proof of the surjectivity of the principal
symbol in \cite{H}, shortly before Definition 18.1.21, shows that there exists $A\in\pscl$\
for which $a$\ is a principal symbol. \cqd
 
\begin{thm}\label{cgc2} If $A\in\psh\cap\acc$\ then every principal symbol $a$\ of $A$\ 
satisfies:
\[
\lim_{t\to+\infty}a(tv)\ \ \mbox{{\em exists for all}}\ \ v\in\sx.
\]
\end{thm}
\pf 
Let us show that $\lim_{t\to+\infty}a(tv)=a_0(v)$, for $a_0=\hat\sigma([A])$\ ($\hat\sigma$\ was
defined in the statement of Corollary~\ref{cgc}). 

Let $t_\alpha>0$\ be a net such that $\lim t_\alpha=+\infty$. 
Since $\m_h$\ is compact, the net $t_\alpha v$\ possesses a convergent subnet. Therefore, we
may suppose that $t_\alpha v$\ converges to $m_0\in\m_s$. We must have $m_0\not\in\tix$, 
because $|t_\alpha v|\to\infty$.

Let $A_n$\ be a sequence in $\pscl$\ converging to $A$, and let $a_n=\hat\sigma([A_n])$. By 
definition of $\hat\sigma$, $A_n$\ has a principal symbol $\tilde a_n$\ satisfying 
$\tilde a_n(tv)=a_n(v)$\ for all $v\in\sx$. By definition of $\sigma$, 
$\sigma(A_n)(m_0)=\lim_\alpha \tilde a_n(t_\alpha v)=a_n(v)$. Since $\sigma$\ and 
$\hat\sigma$\ are continuous, we get: 
\[
\sigma(A)(m_0)=\lim_n\sigma(A_n)(m_0)=\lim_na_n(v)=a_0(v),
\]
for $a_0=\hat\sigma([A])$. By definition of $\sigma$, and because $t_\alpha v\to m_0$,  
\[
\sigma(A)(m_0)=\lim_\alpha a(t_\alpha v),
\]
and, hence, 
\[
\lim_{t\to+\infty}a(tv)=a_0(v),
\]
if $a$\ is a principal symbol for $A$. \cqd

To prove that $\acc$\ is strictly contained in $\ac$\ (using Theorem~\ref{cgc2}\ and
the surjectivity of the principal symbol), it is therefore enough to give an example
of a symbol in $S^0(\tix)$\ such that some of the radial limits in the statement of 
Theorem~\ref{cgc2}\ fail to exist. 

The author learned the following example from Jorge Hounie.

{\bf Example}\ Choose a non-negative $g\in C^\infty(\rn)$, positive on 
$\{\xi;1<|\xi|<2\}$\ and with support contained in $\{\xi;1\leq|\xi|\leq 2\}$. It is 
not difficult to check that
\[
\psi(\xi)=\sum_{k=0}^{\infty}g(2^{-k}\xi).
\]
is a symbol in $S^0(\rn)$\ and that $\lim_{t\to\infty}\psi(t\xi)$\ does not exist for any 
$\xi\neq 0$. If $\varphi$\ is a smooth function with support contained in $\tilde U$, 
for some chart $\chi:U\to\tilde U$, then the pullback to $\tix$\ of 
$\varphi(x)\psi(\xi)$\ by the bundle trivialization induced by $\chi$\ gives an example 
of a symbol $a\in S^0(\tix)$\ for which $\lim_{t\to\infty}a(tv)$\ does not exist
for any $v\in\tix$\ such that $\varphi(\chi(\pi(v)))\neq 0$.

\section{$\acc$\ is a comparison algebra}\label{ca}

Let $\Delta$\ denote the Laplace \op\ on a Riemannian manifold $\Omega$\ equipped with the 
surface measure $d\!\mu$, and let $H$\ denote the {\em Friedrichs extension}\ of $1-\Delta$\ 
(this is a special self-adjoint realization on the Hilbert space $L^2(\Omega,d\!\mu)$\ of the 
symmetric operator $1-\Delta$\ with domain $C_0^{\infty}(\Omega)$). A {\em comparison algebra}\ 
for the triple $\{\Omega,d\!\mu,H\}$\ was defined 
\footnote{His definition is actually more general than this: the measure does not have 
to be the surface measure for some Riemannian metric, the elliptic \op\ does not have to be
$1-\Delta$. He only requires that the \op\ and the measure be related by 
\cite{C2}, (V.1.1).} 
by Cordes (\cite{C2}, Chapter~V)\ as the smallest norm-closed, adjoint-invariant algebra of 
bounded \ops\ on $L^2(\Omega,d\!\mu)$\ containing:
\begin{enumerate}
\item all \ops\ of multiplication by elements of a class of smooth functions, and 

\item all \ops\ of the form $L\Lambda$, where $\Lambda=H^{-1/2}$\ and $L$\ is an element of a 
class of first-order linear partial differential \ops .
\end{enumerate}
These classes of functions and of differential operators must then satisfy certain axioms. 

In case $\Omega=X$\ is compact and without boundary, those axioms imply that there is only one 
comparison algebra: it is the \cst-algebra generated by all multiplications by smooth functions 
and by all \ops\ of the form $L\Lambda$, where $L$\ is any first-order linear differential \op\ 
with smooth coefficients. So, let $\Cc$\ denote the unique comparison algebra over $X$, defined 
by a chosen Riemannian metric. All generators of $\Cc$\ are classical \psd\ \ops; hence
$\Cc\subseteq\acc$. 

Granted that $\Cc$\ contains the compact ideal $\kc$\ (\cite{C2}, Lemma~V.1.1), the equality 
$\acc=\Cc$\ is an easy consequence of following description of the Gelfand space and of the 
Gelfand mapping for the commutative \cst-algebra $\Cc/\kc$. This is a particular case of 
Theorem~VI.2.2 in \cite{C2}.

\begin{thm} The restriction of $\hat\sigma$\ to $\Cc/\kc$\ is
a \cst-algebra isomorphism.
\end{thm}
\pf In view of Corollary~\ref{cgc}, all we have to show is that this restriction
of $\hat\sigma$\ is surjective. It is not difficult to show that the set of all 
functions
\[
a_0(v)=\lim_{t\to\infty}a(tv),
\]
where $a$\ is a principal symbol of an $A=L\Lambda$, $L$\ a smooth vector field,
separates points on $\sx$. The proof of Theorem~\ref{cgc2}\ shows that $a_0=\hat\sigma([A])$.
By the Stone-Weierstrass Theorem, the image of a subalgebra of $\Cc/\kc$\ is therefore dense in 
$C(\sx)$. \cqd

This shows that, for each $A\in\acc$, there is $A^\prime\in\Cc$\ such that $A-A^\prime$\
is compact. Since $\kc\subset\Cc$, we then get $\acc=\Cc$.

\section{Acknowledgements}
The author was deeply influenced by Cordes' views about \psd\ \ops, and about the \cst-algebras 
generated by them, summarized in three books \cite{C1,C2,C3}. In particular, the main idea for 
this paper was borrowed from \cite{C3}, Theorem V.10.3. 

He also thanks Jorge Hounie and Elmar Schrohe for many conversations.

This work was partially supported by the Brazilian agency CNPq (Processo 300330/88-0).

\vskip0.2cm

2000 Mathematics Subject Classification: 35S05 (58J40, 47G30).
\vskip0.2cm
 
Instituto de  Matem\'atica e Estat\'\i stica, Universidade de S\~ao Paulo,

Caixa Postal 66281, 05311-970 S\~ao Paulo, Brazil. 

melo@ime.usp.br.

\end{document}